\documentclass[12pt]{article}
\usepackage[a4paper,margin=1in]{geometry}
\usepackage{amsmath,amssymb,amsthm,mathtools}
\usepackage[T1]{fontenc}
\usepackage{lmodern}
\usepackage{microtype}
\usepackage{xcolor}
\usepackage[colorlinks=true,linkcolor=blue!45!black,citecolor=blue!45!black,urlcolor=blue!45!black]{hyperref}
\hypersetup{pdftitle={The existence of small defining polynomials for number fields, II},
pdfauthor={Manjul Bhargava and Sameera Vemulapalli}}
\newtheorem{theorem}{Theorem}
\newtheorem{lemma}{Lemma}[section]
\newtheorem{proposition}[lemma]{Proposition}

\theoremstyle{remark}

\newcommand{\Q}{\mathbb Q}
\newcommand{\R}{\mathbb R}
\newcommand{\C}{\mathbb C}
\newcommand{\Z}{\mathbb Z}
\newcommand{\OK}{\mathcal O_K}
\newcommand{\N}{\operatorname{N}}
\newcommand{\Tr}{\operatorname{Tr}}
\newcommand{\Disc}{\operatorname{Disc}}
\newcommand{\covol}{\operatorname{covol}}

\newcommand{\Res}{\operatorname{Res}}
\newcommand{\HH}{\mathcal H}
\numberwithin{equation}{section}

\title{\Large The existence of small defining polynomials for number fields, II}
\author{Manjul Bhargava and Sameera Vemulapalli}
\date{\today}

\begin{document}
\maketitle

\begin{abstract}
We prove that every number field of odd degree $n\geq 3$ and
absolute discriminant~$D$ has an integral defining polynomial
whose coefficients are $O_n(D^{(n-1)/(2n-1)})$ in size. This 
strengthens, in odd-degree cases, the recently established bound of $O_n(D^{1/2})$ for all degrees $n$, which was conjectured in 1998 by Ruppert and is now known to be optimal in even degrees. More generally, we prove that 
this strengthened bound holds for all degree-$n$ fields having a real
embedding and no subfield of index two. The proof makes use of the lattice shapes of the ring of integers and its trace dual, 
the inverse different.

Assuming the quasi-Riemann hypothesis for Dedekind zeta
functions,  as recently  announced  by Cureton and Thorner, 
we obtain an exponent smaller than $1/2$ for every field with
no index-two subfield, of any signature. In the case of primitive~fields, 
we obtain the exponent $(n-1)/(2n)$. In particular, we prove that the optimal
exponent for cubic fields is $1/3$.
\end{abstract}

\section{Introduction}

In this article, we prove the following theorem.

\begin{theorem}\label{thm:odd}
Let $K$ be a number field of odd degree $n$ and absolute
discriminant $D_K$. Then there exists a polynomial
\[
 f(x)=a_0x^n+a_1x^{n-1}+\cdots+a_n\in\Z[x]
\]
of degree $n$ with a root $\alpha$ such that $K=\Q(\alpha)$ and
$H(f):=\max\{|a_i|\}\ll_n D_K^{(n-1)/(2n-1)}$.
\end{theorem}

The existence of an integral defining polynomial of height
$O_n(D_K^{1/2})$ for every number field $K$ of degree $n$
was proved in our previous paper \cite{BV} without using QRH,
and independently also by Frei, Sofos, and Widmer
\cite[Corollary~1.7]{FSW} using QRH for the Dedekind zeta
function of the normal closure of $K$. This established a 1998 conjecture of Ruppert \cite[Question~2]{Ruppert}.
For fields with a real embedding---in particular, for fields of odd degree---such a polynomial was already known to exist by the work of Ruppert for totally real fields of prime degree and Vaaler and Widmer \cite{VW} in general. 

Theorem~\ref{thm:odd} gives an exponent below \(1/2\) for every odd degree \(n\geq3\), improving upon the bound in Ruppert’s conjecture.
In contrast, Widmer \cite[Theorem~2]{Widmer} proved that
the exponent $1/2$ cannot be decreased for any even degree.
The obstruction occurs along certain sequences of quadratic extensions of a fixed number~field.
Indeed, for a fixed number field $F$ and 
quadratic extensions $K/F$, Silverman's relative discriminant
inequality~\cite{Silverman} gives $H(f)\gg_F D_K^{1/2}$
for every defining polynomial $f$ of $K$; see
\cite[Proposition~1]{Widmer}.

We also prove a number of additional results using similar methods. 
First, the
bound in Theorem~\ref{thm:odd} holds more generally for number fields
$K$ having a real embedding and no subfield of index two.
Let us write 
\[
 \HH(K):=\min\{H(f):f\in\Z[x],\ \deg f=n,\
                         K\cong\Q[x]/(f)\}.
\]
Then we prove the following theorem. 
\begin{theorem}\label{thm:real}
Let $K$ be a number field of degree $n$ having a real
embedding and no subfield of index~two. Then
\[
 \HH(K)\ll_n D_K^{(n-1)/(2n-1)}.
\]
\end{theorem}
\noindent Theorem~\ref{thm:real} implies Theorem~\ref{thm:odd}. 
 
For totally real fields, we are able to improve the bound in
Theorem~\ref{thm:real} even further.
For $K\ne\Q$, let 
\[
 m(K):=\max_{\Q\subseteq F\subsetneq K}[F:\Q].
\]
Thus $m(K)=1$ precisely when $K$ is {\em primitive}, that is, when
it has no proper subfield other than $\Q$.
Write
\[
 \rho(n,m):=\frac{n-1}{2n-1+1/(2m-1)}.
\]

\begin{theorem}\label{thm:totally-real}
Let $K$ be a totally real number field of degree $n\geq 3$
with no subfield of index~two. Then
\[
 \HH(K)\ll_n D_K^{\rho(n,m(K))}.
\]
In particular, if $K$ is primitive, then
\[
 \HH(K)\ll_n D_K^{(n-1)/(2n)}.
\]
\end{theorem}

For every $m\geq1$, we have
$\rho(n,m)<(n-1)/(2n-1)$, so Theorem~\ref{thm:totally-real} gives an improvement over Theorem~\ref{thm:real} in the case of totally real fields.  We note that the implied constants in Theorems~\ref{thm:odd}, \ref{thm:real}, and~\ref{thm:totally-real} are all effectively computable from our proofs. 

Our proofs use two complementary ways to construct a
defining polynomial. A short integral generator produces a
monic norm polynomial. Short elements $\alpha,\beta$ in the
inverse different instead produce the integral polynomial
$D_K\cdot\N_{K/\Q}(\alpha x-\beta)$, which need not be monic.
The lattices $\OK$ and $\mathfrak d_K^{-1}$ are dual under
the trace pairing. Inequalities between their successive
minima show that if $\OK$ does not have suitable 
short elements, then the inverse different does; 
the absence of an index-two subfield ensures that enough of
these elements are available to form a ratio generating $K$.
Keeping track of their lengths gives complementary height
bounds. An integer translation at a real embedding saves
a factor in the integral construction, and balancing the
two bounds then proves Theorems~\ref{thm:odd} and~\ref{thm:real}.

For the stronger bounds in the case of totally real fields, 
we distinguish the length of a
shortest element of the inverse different from the length
needed to find a generating ratio. The trace-zero sublattice
enables us to prove an additional relation between these lengths. When
$K$ is totally real, translation at a suitable embedding in conjunction with this additional relation allows us to prove Theorem~\ref{thm:totally-real}.

If we assume the recent results towards the quasi-Riemann
hypothesis (QRH) for Dedekind zeta functions, we can improve the exponents in
Theorems~\ref{thm:odd} and~\ref{thm:real} further, and
extend the stronger bound in Theorem~\ref{thm:totally-real}
to fields of every signature. Following the recent work of
OpenAI \cite{OpenAI} on QRH for the Riemann zeta function, Cureton and Thorner
\cite[Theorem~1.1]{CT} prove that the Dedekind zeta function
of every number field has no zeros in $\Re(s)>31/32$.
In the theorems that follow, we state the QRH 
hypothesis explicitly with a parameter
$N\geq2$, so that the result of Cureton and Thorner
corresponds to $N=32$ and the dependence of our constants
is clear.

We begin with the following theorem that strengthens Theorems~\ref{thm:odd}--\ref{thm:totally-real}. 

\begin{theorem}\label{thm:qrh}
Fix $N\geq2$. Let $K$ be a number field of degree $n\geq 3$
with no subfield of index~two. If its Dedekind zeta function
$\zeta_K(s)$ has no zeros in the half-plane
$\Re(s)>1-1/N$, then
\[
 \HH(K)\ll_{n,N}D_K^{\rho(n,m(K))}.
\]
In particular, if $K$ is primitive, then
\[
 \HH(K)\ll_{n,N}D_K^{(n-1)/(2n)}.
\]
\end{theorem}

Only the implied constant depends on $N$; the exponent is independent of $N$. The hypothesis concerns
$\zeta_K$ itself; no assumption about the zeta function of the 
normal closure of $K$ is needed. The implied constant is effectively computable
in terms of $n$ and $N$.

Finally, in degree 3, by combining Theorem~\ref{thm:qrh} with the work of Dubickas~\cite{Dubickas} on cubic lower bounds, we can show that the exponent in Theorem~\ref{thm:qrh} is optimal. 

\begin{theorem}\label{thm:cubic}
Fix $N\geq2$, and suppose that the Dedekind zeta functions
of cubic fields have no zeros in $\Re(s)>1-1/N$. Then
\[
 \HH(K)\ll_N D_K^{1/3}
\]
for every cubic field $K$. Moreover, there are infinitely
many cubic fields $K$ for which
\[
 \HH(K)\geq\frac16D_K^{1/3}.
\]
Hence $1/3$ is the optimal exponent for cubic fields $K$. 
\end{theorem}

The proof of Theorem~\ref{thm:qrh} (and Theorem~\ref{thm:cubic}) initially uses the same two norm
constructions. The zero-free half-plane is then used to supply
a degree-one prime $p$ at a suitable scale. A root $r$ modulo
$p$ allows us to replace $f(x)$ by $f(px+r)/p$,
saving a factor of~$p$ while preserving the field. Applying
this operation to both norm constructions proves
Theorem~\ref{thm:qrh}.

We also give explicit lower-bound constructions. The cubic
lower exponent is due to Dubickas \cite[Theorem~1]{Dubickas}.
In a suitable
pure cubic field, either a generator or its reciprocal has
a nonzero quadratic coefficient in the standard integral
basis. Integrality then forces a large conjugate and hence
a large Mahler measure. A relative version gives families
with $\HH(K)\asymp D_K^{1/3}$ in every degree divisible by
three. 

Theorems~\ref{thm:odd}--\ref{thm:qrh} show that the exponent $1/2$, as in Ruppert's Conjecture~\cite{Ruppert}, 
is not optimal for fields of odd degree, and indeed for fields not having a subfield of index two. 
This raises the question as to what the optimal exponent is in each degree.  For even degrees $n$,
and for $n=3$, the answer is 1/2 and 1/3, respectively, by \cite{BV,FSW} and Theorem~\ref{thm:cubic}.  
More generally, suppose we fix the field degree $n$, the associated Galois group $G\subset S_n$, and the real signature $(r_1,r_2)$. What is the optimal exponent over all fields $K$ having this degree, Galois group, and real signature?  Again, Theorems~\ref{thm:odd}--\ref{thm:qrh} demonstrate that the answer is strictly less than 1/2 in most cases. 

For more on the history of this problem, see \cite{BV}.

This paper is organized as follows. Section~\ref{sec:first}
recalls the inverse different construction and sketches the
proof of the square-root bound in \cite{BV}.
Section~\ref{sec:odd} proves Theorems~\ref{thm:odd},
\ref{thm:real}, and~\ref{thm:totally-real} using successive minima in 
the ring of integers and the inverse different.
These two sections do not make use of any quasi-Riemann hypotheses. Section~\ref{sec:qrh}
then uses the zero-free hypothesis of Cureton and Thorner~\cite{CT} to prove
Theorem~\ref{thm:qrh}. Finally, in Section~\ref{sec:cubic},
we prove the cubic lower bound, deduce Theorem~\ref{thm:cubic},
and carry out the aforementioned relative lower-bound construction.

\section{The inverse different and the square-root bound}
\label{sec:first}

We begin with the constructions from \cite{BV} that are 
used throughout. The case $n=1$ is immediate, so we henceforth
suppose that $K$ has degree $n\geq2$, and write $D=D_K$.
As in \cite{BV}, choose all $r_1$ real embeddings and one
embedding from each of the $r_2$ complex conjugate pairs,
and identify $K\otimes_\Q\R$ with
$\R^{r_1}\times\C^{r_2}$. We use the usual Euclidean norm
on this space, identifying each complex coordinate with
$\R^2$, and write $\|\alpha\|$ for the norm of the image
of $\alpha$. Thus
\[
 \|\alpha\|^2=\sum_{v\text{ real}}|\sigma_v(\alpha)|^2
                 +\sum_{v\text{ complex}}|\sigma_v(\alpha)|^2.
\]
We write $\N(I)$ for the absolute norm of a fractional
ideal, and $\N(\alpha)=\N_{K/\Q}(\alpha)$ for the field
norm of an element. Then the discriminant formula gives
\begin{equation}\label{eq:covol}
 \covol(I)=2^{-r_2}\sqrt D\,\N(I);
\end{equation}
see \cite[Chapter~2]{Marcus}. Each of the $n$ conjugates
of $\alpha$, including both members of a complex pair,
has absolute value at most $\|\alpha\|$.

\begin{lemma}\label{lem:norm}
Let $I$ be a nonzero fractional ideal of $K$, and let
$\alpha,\beta\in I$. Then
\[
 F(x,y)=\frac{\N_{K/\Q}(\alpha x+\beta y)}{\N(I)}
 \in\Z[x,y].
\]
The coefficient of $x^iy^{n-i}$ is at most
\[
 \binom ni\frac{\|\alpha\|^i\|\beta\|^{n-i}}{\N(I)}
\]
in absolute value. Moreover,
\begin{equation}\label{eq:product-height}
 H(F)\leq\frac{2^n}{\N(I)}
 \prod_{\sigma:K\hookrightarrow\C}
 \max\{|\sigma(\alpha)|,|\sigma(\beta)|\}.
\end{equation}
If $\alpha\ne0$ and $\beta/\alpha$ generates $K$, then
$F(x,-1)$ is a defining polynomial for $K$.
\end{lemma}

\begin{proof}
Multiplication by $\alpha x+\beta y$ is a map from
$\OK[x,y]$ to $I[x,y]$. In integral bases its matrix has
entries in $\Z[x,y]$, and its determinant is
$\N(\alpha x+\beta y)/\N(I)$ up to sign. This proves
integrality. Expanding the product over the $n$ embeddings
gives the coefficient bounds. The final assertion follows
because $F(\beta/\alpha,-1)=0$ and $F(x,-1)$ has degree~$n$.
\end{proof}

The following lemma shows that a bounded
combination of short elements gives a generating ratio.

\begin{lemma}\label{lem:generate}
Let $V\subseteq K$ be a $k$-dimensional $\Q$-vector space
with basis $\alpha_1,\ldots,\alpha_k$, and suppose that
$k>m(K)$. For every $0\ne\alpha\in V$, there are integers
$0\leq c_i\leq\binom n2$ such that
\[
 \beta=\sum_{i=1}^k c_i\alpha_i,
 \qquad \Q(\beta/\alpha)=K.
\]
In particular, if $1\in V$ and $V$ has a basis in $\OK$
of length at most $X$, then $K$ has an integral generator
of length $O_n(X)$.
\end{lemma}

\begin{proof}
The ratios $\alpha_i/\alpha$ span a space of dimension $k$.
The field they generate cannot be proper, since every proper
subfield has degree at most $m(K)<k$. For distinct embeddings
$\sigma,\tau$ of $K$, the equality
\[
 \sum_i c_i\bigl(\sigma(\alpha_i/\alpha)
                   -\tau(\alpha_i/\alpha)\bigr)=0
\]
is therefore a nonzero linear equation in the $c_i$. Put
$M=\binom n2$. Each equation contains at most $(M+1)^{k-1}$
points of $\{0,\ldots,M\}^k$, and there are only $M$ such
equations. Their union does not exhaust this set. A point
outside the union gives a ratio with $n$ distinct conjugates,
and hence a generator. The last assertion follows by taking
$\alpha=1$.
\end{proof}

Let $\mathfrak d_K$ be the different, and put
$I^\vee=\mathfrak d_K^{-1}I^{-1}$. The trace pairing
identifies $I^\vee$ with the integral dual of $I$; see
\cite[Chapter~3, Exercises~33--39]{Marcus} for these
identities and the trace-dual basis construction.
With the Euclidean norm just defined, the Euclidean dual of
$I$ is obtained from $I^\vee$ by conjugating each complex
coordinate and multiplying it by two. Indeed, the trace of
a product is the sum of its real coordinates and twice the
real parts of its complex coordinates. Conjugation is an
isometry, and the extra factors lie between one and two.
Consequently, the successive minima of the trace dual and
the Euclidean dual differ by at most a factor of two. These
maps act on the ambient real vector space; conjugation here
need not be an automorphism of $K$.

Write $\lambda_j(L)$ for the least radius of a ball about
the origin containing $j$ linearly independent vectors of
a lattice $L$. Let $L^*$ be its Euclidean dual, namely, the
lattice of vectors having integral inner product with every
vector of $L$. We shall repeatedly use the inequalities
\begin{equation}\label{eq:dual-minima}
 1\leq\lambda_j(L)\lambda_{n+1-j}(L^*)\ll_n1.
\end{equation}
These follow directly from Minkowski's second theorem;
see \cite[Chapter~VIII, \S4.3, Theorem~V]{Cassels}.
We give the argument to
make clear how the two lattices are related. Choose $j$
independent vectors in $L$ of length at most $\lambda_j(L)$
and $n+1-j$ independent vectors in $L^*$ of length at most
$\lambda_{n+1-j}(L^*)$. If the product of these two lengths
were less than one, all their inner products would be integers
of absolute value less than one, and hence zero. Their spans
would then be orthogonal, which is impossible because their
dimensions sum to $n+1$. This proves the lower bound.

For the upper bound, let $\omega_n$ be the volume of the
Euclidean unit ball. Minkowski's second theorem gives
\[
 \prod_{i=1}^n\lambda_i(L)
 \leq\frac{2^n}{\omega_n}\covol(L).
\]
Apply the same inequality to $L^*$ and use
$\covol(L^*)=\covol(L)^{-1}$. The product of the $n$
paired quantities in \eqref{eq:dual-minima} is at most
$(2^n/\omega_n)^2$. Since each of them is at least one,
each is bounded above by the same constant. Thus a large
successive minimum in one lattice forces the corresponding
minimum in the dual lattice to be small.

We also use the elementary effective bound
\begin{equation}\label{eq:elementary-height}
 \HH(K)\ll_n D^{n/2}.
\end{equation}
Indeed, a nonzero algebraic integer $\gamma$ satisfies
$1\leq|\N(\gamma)|\leq\|\gamma\|^n$, so every successive
minimum of $\OK$ is at least one. Minkowski's second theorem
therefore gives $\lambda_n(\OK)\ll_n\sqrt D$.
Lemma~\ref{lem:generate}, applied to $V=K$ and $\alpha=1$,
provides an integral generator of length $O_n(\sqrt D)$.
Its monic norm polynomial has the required height. This
bound will suffice whenever we need to absorb a bounded
range of discriminants.

We now recall how these constructions give the square-root
bound of \cite{BV}. We normalize an ideal lattice to covolume
one, allowing positive rescalings of the archimedean
coordinates with product one over the $n$ embeddings.
Its Euclidean dual is, up to conjugation of the complex
coordinates, another lattice in this normalized family. For
$k=\lceil n/2\rceil$, \eqref{eq:dual-minima} implies that
one of these lattices has $k$ independent vectors of bounded
length. If $K$ has no index-two subfield, then
$m(K)\leq n/3<k$. Lemmas~\ref{lem:norm} and
\ref{lem:generate} give a defining polynomial of height
$O_n(\sqrt D)$: the archimedean rescalings cancel in the
product in \eqref{eq:product-height}, and the normalization
contributes a factor $2^{-r_2}\sqrt D$.

The remaining argument in \cite{BV} uses subconvexity.
If the construction fails for every ideal and every
archimedean rescaling, the short elements of each ideal
span a line over an index-two subfield $F$. For a fixed
ideal this line is independent of the rescaling and is
preserved by all units of $K$. Dirichlet's unit theorem
forces $F$ to be totally real and $K/F$ to be CM.
Writing $n=2s$ and $D=D_F^2d$, a count of relative ideal
classes and the class-number formula then give
\begin{equation}\label{eq:small-base}
 D_F\ll_{s,\epsilon}D^\epsilon\qquad(\epsilon>0);
\end{equation}
see \cite[Sections~3--4]{BV}.

For these CM fields, let $L$ be a normalized ideal lattice
and define the completed Epstein zeta function by
\[
 Z_L(z)=\pi^{-nz/2}\Gamma(nz/2)
                \sum_{0\ne v\in L}\|v\|^{-nz}
 \qquad(\Re(z)>1),
\]
with meromorphic continuation elsewhere. Poisson summation
gives a constant $r_s>0$ such that $Z_L(1/2)\geq1$
whenever $\lambda_1(L)\leq r_s$. We average over ideal classes
with equal weights and over volume-preserving rescalings
modulo units with probability measure. Hecke's formula
\cite{Hecke,Khayutin} gives this average as
\[
 c_sD^{-1/4}\frac{\zeta_K(1/2)}{\kappa_K},
 \qquad \kappa_K=\Res_{z=1}\zeta_K(z),\qquad c_s>0.
\]
The subconvexity theorem of Michel and Venkatesh \cite{MV}
gives $|\zeta_K(1/2)|\ll_s D_F^a d^{1/4-\eta}$ for some
$a\geq0$ and $\eta>0$ depending only on $s$. Together with
\eqref{eq:small-base} and the residue estimate
$\kappa_K^{-1}\ll_{n,\epsilon}D^\epsilon$ recalled in
\cite[Section~4]{BV}, this makes the average tend to zero.
Some lattice therefore has $\lambda_1(L)>r_s$; Minkowski's
second theorem bounds all its successive minima, yielding
the required generating ratio. This completes the sketch
of the square-root bound. The new results below do not
depend on this analytic argument.

\section{Bounds from successive minima}
\label{sec:odd}

Let $m=m(K)$, so a
space of dimension $m+1$ cannot lie in a proper subfield.
The $(m+1)$st minimum of $\OK$ measures how far we must go
to find enough independent integers to construct a generator.
If this minimum is small, the integers themselves give a
small polynomial. If it is large, \eqref{eq:dual-minima}
provides short elements of the inverse different, and these
give a different small polynomial. A real embedding improves
the first construction by allowing one conjugate to be made
bounded by an integer translation.

\begin{proof}[Proof of Theorems~$\ref{thm:odd}$ and~$\ref{thm:real}$]
We first prove Theorem~\ref{thm:real}. Let $K$ have a real
embedding and no index-two subfield. Put $m=m(K)$, $k=m+1$, and
\[
 X=\lambda_k(\OK).
\]
Every proper subfield has index at
least three. Thus $m\leq n/3$; in particular, $2m<n$ and
$k\leq n+1-k$. Every nonzero algebraic integer has length
at least one,
as observed above, so $X\geq1$. Choose $k$ independent
elements of $\OK$ of length at most $X$. We may arrange
for their span to contain $1$: if it does not, replace one
of them by $1$. The resulting vectors remain independent
and have length $O_n(X)$, since $\|1\|\leq\sqrt n$.
Lemma~\ref{lem:generate}, with $\alpha=1$, then gives an
integral generator $\theta$ of length $O_n(X)$.

We choose a real embedding $\sigma_0$ of $K$, and let $a\in\Z$
be a nearest integer to $\sigma_0(\theta)$. Then
\[
 |\sigma_0(\theta-a)|\leq\tfrac12,
 \qquad |\sigma(\theta-a)|\ll_n X
\]
for every other embedding $\sigma$. Indeed,
$|a|\leq|\sigma_0(\theta)|+1/2\ll_n X$, which gives the
second estimate. Translation by $a$ preserves both integrality
and the generated field. The product bound in
Lemma~\ref{lem:norm}, applied to $I=\OK$, $\alpha=1$,
and $\beta=\theta-a$, now yields
\[
 H\bigl(\N_{K/\Q}(x-(\theta-a))\bigr)
 \leq 2^n\prod_\sigma\max\{1,|\sigma(\theta-a)|\}
 \ll_n X^{n-1}.
\]
There are only $n-1$ potentially large factors, since the
factor at $\sigma_0$ is one. We have proved
\begin{equation}\label{eq:primal-real}
 \HH(K)\ll_n X^{n-1}.
\end{equation}

For the complementary construction, take
$I=\mathfrak d_K^{-1}$. Its successive minima are within
a factor of two of those
of the Euclidean dual of $\OK$. Equation~\eqref{eq:dual-minima}
and $n+1-k\geq k$ therefore give
\[
 \lambda_k(I)
 \ll_n\lambda_{n+1-k}(\OK)^{-1}
 \leq X^{-1}.
\]
Choose $k$ independent elements of $I$ of length $O_n(X^{-1})$
and take any one of them as $\alpha$. Their span has dimension
$k>m$, so Lemma~\ref{lem:generate} supplies a bounded
integer combination $\beta$ for which $\beta/\alpha$
generates $K$. Both elements have length $O_n(X^{-1})$.
Since $\N(I)=D^{-1}$, Lemma~\ref{lem:norm} shows that
$D\N(\alpha x-\beta)$ is an integral defining polynomial
whose coefficients are $O_n(DX^{-n})$. Consequently,
\begin{equation}\label{eq:dual-basic}
 \HH(K)\ll_n D X^{-n}.
\end{equation}

Combining \eqref{eq:primal-real} and \eqref{eq:dual-basic},
we obtain
\[
 \HH(K)\ll_n\min\{X^{n-1},DX^{-n}\}
 \leq D^{(n-1)/(2n-1)}.
\]
Indeed, use the first bound when $X\leq D^{1/(2n-1)}$
and the second otherwise. This proves Theorem~\ref{thm:real}.
Finally, the complex embeddings
of a number field occur in conjugate pairs. An odd-degree
field therefore has a real embedding; the tower law also
precludes an index-two subfield. This proves Theorem~\ref{thm:odd}.
\end{proof}

The proof shows exactly where each hypothesis enters.
The subfield condition ensures that sufficiently many short
elements have a generating ratio. The real embedding saves
one factor of $X$ in \eqref{eq:primal-real}. Integer
translation at a complex embedding does not control its
imaginary part, so that step does not give the same bound
for a totally complex field.

For the inverse different construction, we distinguish
the length of a shortest nonzero element from the length
needed to find a generating ratio. The additional information
comes from the trace-zero sublattice, whose covolume is smaller
than that of $\OK$ by a factor of order $D$.

\begin{lemma}\label{lem:tracezero}
Let $I=\mathfrak d_K^{-1}$ and
\[
 \Lambda=\{\alpha\in I:\Tr_{K/\Q}(\alpha)=0\}.
\]
Then $\Lambda$ has rank $n-1$ and
\begin{equation}\label{eq:trace-covol}
 \covol(\Lambda)=2^{-r_2}\sqrt{r_1+4r_2}\,D^{-1/2}.
\end{equation}
Let $m=m(K)<n/2$, $k=m+1$, and $X=\lambda_k(\OK)$.
If $X$ exceeds a sufficiently large constant depending only
on $n$, and
\[
 A=\lambda_1(\Lambda),\qquad B=\lambda_k(\Lambda),
\]
then
\begin{equation}\label{eq:AB}
 B\ll_n X^{-1},\qquad
 A^mB^{n-m-1}\ll_n D^{-1/2},\qquad
 A\geq D^{-1/n}.
\end{equation}
There are $\alpha,\beta\in\Lambda$ such that
\[
 \|\alpha\|=A,\qquad \|\beta\|\ll_n B,
 \qquad \Q(\beta/\alpha)=K.
\]
\end{lemma}

\begin{proof}
First, $1$ is primitive as a vector of $\OK$: if $1=q\gamma$
with $q\in\Z_{>0}$ and $\gamma\in\OK$, then $1/q$ is a
rational algebraic integer, so $q=1$. We may therefore extend
$1$ to an integral basis of $\OK$. Its trace-dual basis is
an integral basis of $I$ and contains an element $\eta$ with
$\Tr(\eta)=1$. Since every element of $I$ has integral
trace, this proves $\Tr(I)=\Z$ and
\[
 I=\Lambda\oplus\Z\eta.
\]
In particular, $\Lambda$ has rank $n-1$.

In our Euclidean coordinates, the trace is the inner
product with the vector $\mathbf t$ whose real coordinates
are $1$ and whose complex coordinates are $2$.
Thus $\|\mathbf t\|=\sqrt{r_1+4r_2}$, and the trace-zero
hyperplane is perpendicular to $\mathbf t$. The perpendicular
component of $\eta$ has length
$\Tr(\eta)/\|\mathbf t\|=1/\sqrt{r_1+4r_2}$.
The volume of a fundamental parallelepiped is its base
volume times this height, so
\[
 \covol(I)=\frac{\covol(\Lambda)}{\sqrt{r_1+4r_2}}.
\]
Since \eqref{eq:covol} gives
$\covol(I)=2^{-r_2}D^{-1/2}$, this proves
\eqref{eq:trace-covol}.

Next, \eqref{eq:dual-minima} gives $k$ independent elements
of $I$ of length $O_n(X^{-1})$, just as in
the proof above. For each such element $\gamma$,
Cauchy--Schwarz gives
$|\Tr(\gamma)|\leq\sqrt{r_1+4r_2}\,\|\gamma\|$.
When $X$ exceeds a sufficiently large constant depending
only on $n$, this is less than one. Its integral trace must
therefore vanish, so all $k$ elements belong to $\Lambda$.
This proves $B\ll_n X^{-1}$.

Let $\mu_1\leq\cdots\leq\mu_{n-1}$ be the successive
minima of $\Lambda$. The first $m$ are at least
$A=\mu_1$, while each of the remaining $n-m-1$ is at least
$B=\mu_{m+1}$. Minkowski's second theorem, now in dimension
$n-1$, thus gives
\[
 A^mB^{n-m-1}
 \leq\prod_{j=1}^{n-1}\mu_j
 \ll_n\covol(\Lambda)
 \ll_n D^{-1/2}.
\]
For $0\ne\alpha\in I$, Lemma~\ref{lem:norm} shows that
$D\N(\alpha)$ is a nonzero integer. Since
$|\N(\alpha)|\leq\|\alpha\|^n$, we have
$\|\alpha\|\geq D^{-1/n}$.

Finally, choose a shortest nonzero vector $\alpha$ of
$\Lambda$. Since the vectors in the ball of radius $B$
span a space of dimension at least $k$, we can extend
$\alpha$ to $k$ independent vectors in that ball.
Their span has dimension $k=m+1>m(K)$, so
Lemma~\ref{lem:generate} supplies $\beta$ as a bounded
integer combination of these vectors. It lies in $\Lambda$,
has length $O_n(B)$, and its ratio to $\alpha$ generates $K$.
\end{proof}

The estimates in this lemma will be used in both of the
remaining upper bounds. Put
\begin{equation}\label{eq:exponent-parameters}
 a=\frac{2m-1}{2m},\qquad
 b=\frac{n(m-1)+1}{m}.
\end{equation}
Solving the product bound in \eqref{eq:AB} for $A$, and
then using $B\ll_n X^{-1}$, gives
\begin{align}
 DA B^{n-1}
 &\ll_n D^{(2m-1)/(2m)}B^{[n(m-1)+1]/m}\notag\\
 &\ll_n D^aX^{-b}.\label{eq:dual-refined}
\end{align}
The exponent of $B$ is positive. Moreover,
\begin{equation}\label{eq:balance}
 \min\{X^{n-1},D^aX^{-b}\}
 \leq D^{a(n-1)/(n-1+b)}=D^{\rho(n,m)}.
\end{equation}
To see the inequality, use the first term when
$X\leq D^{a/(n-1+b)}$ and the second otherwise.
The identity of the exponents follows from
\[
 \frac{a(n-1)}{n-1+b}
 =\frac{(2m-1)(n-1)}{2\bigl(n(2m-1)-m+1\bigr)}
 =\frac{n-1}{2n-1+1/(2m-1)}.
\]

\begin{proof}[Proof of Theorem~$\ref{thm:totally-real}$]
Put $m=m(K)$, $k=m+1$, and $X=\lambda_k(\OK)$.
If $X$ is bounded, the short integral generator already
gives bounded height. Otherwise Lemma~\ref{lem:tracezero}
provides $\alpha,\beta,A,B$.
The integral construction already gives
$\HH(K)\ll_n X^{n-1}$ by \eqref{eq:primal-real}.
For the inverse different bound,
choose an embedding $\sigma_0$ with
$|\sigma_0(\alpha)|\geq A/\sqrt n$; such an embedding
exists because $\sum_\sigma|\sigma(\alpha)|^2=A^2$.
The totally real condition ensures that this particular embedding is real.
We may therefore choose an integer $t$ nearest to
$\sigma_0(\beta)/\sigma_0(\alpha)$. Since $B\geq A$, we have 
\[
 |t|\leq\frac{|\sigma_0(\beta)|}{|\sigma_0(\alpha)|}
          +\frac12\ll_n B/A,
\]
and $\beta'=\beta-t\alpha$ satisfies
\[
 |\sigma_0(\beta')|\ll_n A,
 \qquad |\sigma(\beta')|\ll_n B
 \quad\text{for all other }\sigma.
\]
The ratio $\beta'/\alpha=\beta/\alpha-t$ still generates
$K$. At $\sigma_0$, both $\alpha$ and $\beta'$ have
absolute value $O_n(A)$; at every other embedding they have
absolute value $O_n(B)$. Equation~\eqref{eq:product-height}
therefore gives
\[
 H\bigl(D\N(\alpha x-\beta')\bigr)
 \ll_n DA B^{n-1}.
\]
Combining this bound with \eqref{eq:dual-refined} and
\eqref{eq:balance} proves the theorem. \end{proof}

\section{An improvement from a zero-free half-plane}\label{sec:qrh}

An integer translation makes one real conjugate small. To
obtain a comparable saving without a real embedding, we use
a congruence instead: if a defining polynomial $f$ has a root
$r$ modulo a prime $p$, then $f(px+r)/p$ is again an integral
defining polynomial. The division by $p$ is useful provided
the substitution does not increase the other coefficients
too much. We shall choose $p$ to match the lengths of the
elements occurring in the norm polynomial.

The prime-ideal estimate needed for this choice is a direct
consequence of the prime ideal theorem under a fixed
zero-free hypothesis.
We use it only at scales bounded above and below by fixed
positive powers of the discriminant.

\begin{lemma}\label{lem:prime}
Fix $n\geq3$, $N\geq2$, $\epsilon>0$, and $C>0$.
Suppose that $[K:\Q]=n$ and $\zeta_K(s)$ has no zeros
in $\Re(s)>1-1/N$. For $D$ sufficiently large in terms of
these constants, every interval $(Y,2Y]$ with
\[
 D^\epsilon\leq Y\leq D^C
\]
contains $\gg_{n,N,\epsilon,C}Y/\log Y$ rational primes
having a degree-one prime in $K$.
If $f\in\Z[x]$ defines $K$ and $H(f)\ll_n D^C$, one
may choose such a prime not dividing the leading coefficient
or discriminant of $f$. For this choice, $f$ has a root
modulo $p$.
\end{lemma}

\begin{proof}
We apply \cite[Corollary~1.2]{FSW} to the trivial extension
$K/K$, with its parameter equal to $2Y$. The zero-free
width in that result is $1/N$,
and its discriminant threshold is
$O_N((\log D)^N(\log\log D)^{2N})$. This is smaller
than $D^\epsilon$ for sufficiently large $D$, with an
effectively computable threshold. It gives $\gg Y/\log Y$
prime ideals with norms
in $(Y,2Y]$. Only $\zeta_K$ occurs in this application;
there is no passage to a normal closure.
Prime ideals of residue degree at least two lie above
rational primes at most $\sqrt{2Y}$, so there are only
$O_n(\sqrt Y)$ of them. A rational prime has at most $n$
prime ideals above it, which proves the first assertion.

Write $a$ for the leading coefficient of $f$. Its discriminant
is nonzero, since $f$ is irreducible in characteristic zero,
and is a homogeneous polynomial of degree $2n-2$ in the
coefficients. Thus $|a\Disc(f)|=D^{O_{n,C}(1)}$, and this
integer has $O_{n,C}(\log D)$ distinct prime divisors.
Since $Y\geq D^\epsilon$, deleting them does not exhaust
the available primes.

We explain the last assertion also when $f$ is nonmonic.
Let $\theta$ be a root of $f$ and let $p\nmid a\Disc(f)$.
The polynomial $a^{-1}f$ is monic over $\Z_p$, so
$\theta$ is integral over $\Z_p$. Moreover, the
discriminant of the basis $1,\theta,\ldots,\theta^{n-1}$
is a $p$-adic unit. The discriminant-index formula therefore
shows that
\[
 \Z_p[\theta]=\OK\otimes_\Z\Z_p.
\]
Reducing this identity modulo $p$, the factorization of
$a^{-1}f$ has the residue degrees of the primes of $K$
above $p$. This is Dedekind's factorization theorem; see
\cite[Theorem~4.8.13]{Cohen}. A prime of residue degree one
therefore supplies a linear factor, and hence a root of $f$
modulo $p$.
\end{proof}

If $r$ is a root of $f$ modulo $p$, with $0\leq r<p$, then
\begin{equation}\label{eq:rescale}
 g(x)=\frac{f(px+r)}p\in\Z[x]
\end{equation}
defines the same field as $f$. Indeed, its constant term is
integral because $p\mid f(r)$. Every contribution to the
coefficient of $x^j$, for $j\geq1$, contains a factor $p^j$,
so division by $p$ leaves all other coefficients integral
as well. The leading coefficient remains nonzero. If
$\theta$ is a root of $f$, then $(\theta-r)/p$ is a root
of $g$ and generates the same field. Since $g$ still has
degree $n$, it is a defining polynomial for $K$.

\begin{proof}[Proof of Theorem~\ref{thm:qrh}]
Put $m=m(K)$, $k=m+1$, $X=\lambda_k(\OK)$, and
$e=\rho(n,m)$. Let $a,b$ be as in
\eqref{eq:exponent-parameters}. We may suppose that $D$
exceeds a sufficiently large effective constant depending
only on $n$ and $N$: the elementary bound
\eqref{eq:elementary-height} absorbs the remaining bounded
range into the implied constant. If $X$ is bounded in terms
of $n$, Lemma~\ref{lem:generate}, after including $1$ among
the short vectors as before, gives a defining polynomial
of bounded height. We may therefore also assume that
Lemma~\ref{lem:tracezero} applies.

Choose a positive $\epsilon$, depending only on $n$ and
$m$, sufficiently small that
\begin{equation}\label{eq:epsilon}
 n\epsilon<e,
 \qquad
 \frac{n-2}{2(n-1)}+\frac{nm}{n-1}\epsilon<e.
\end{equation}
This is possible because
\[
 e-\frac{n-2}{2(n-1)}
 =\frac{n(m-1)+1}
 {2(n-1)\bigl(n(2m-1)-m+1\bigr)}>0.
\]
For fixed $n$ there are only finitely many possible values
of $m$. We may therefore make these choices uniformly in
$m$, so that all final constants depend only on $n$ and
$N$. The purpose of $\epsilon$ is to separate the
scales at which the prime lemma applies from the small scales
where the desired height bound already follows directly.

We first obtain the integral bound. As before, there is an
integral generator $\theta$ with $\|\theta\|\ll_n X$.
If $X<D^\epsilon$, its monic minimal polynomial has height
$O_n(X^n)=O_n(D^e)$, and we are done. Otherwise, we apply
Lemma~\ref{lem:prime} to
\[
 f_0(x)=\N_{K/\Q}(x-\theta)
\]
with $Y=X$, giving $X<p\leq2X$; the lemma applies 
because all successive minima of $\OK$ are at least
one, while their product is $O_n(\sqrt D)$ by
Minkowski's second theorem. In particular,
$X\ll_n\sqrt D$ and $H(f_0)\ll_n X^n$, so both the
scale and the polynomial height are bounded by a fixed
power of $D$. If $0\leq r<p$ is a root modulo $p$, each
factor $px+r-\sigma(\theta)$ has coefficients of size
$O_n(p+X)$. Expanding their product and dividing by $p$
gives
\[
 H\left(\frac{f_0(px+r)}p\right)
 \ll_n\frac{(p+X)^n}{p}\ll_n X^{n-1}.
\]
Thus, unless the desired conclusion already holds,
\begin{equation}\label{eq:primal-qrh}
 \HH(K)\ll_{n,N} X^{n-1}.
\end{equation}

We now take $\alpha,\beta,A,B$ as in
Lemma~\ref{lem:tracezero}, and set
\[
 R=B/A\geq1,\qquad
 f_1(x)=D\N_{K/\Q}(\alpha x-\beta).
\]
This polynomial is integral and defines $K$, and its height
is $O_n(DB^n)$. The natural prime scale is now $R=B/A$:
if $p\asymp R$, the elements $p\alpha$ and $\beta$ have
comparable lengths. The substitution $x\mapsto px+r$ then
keeps the factors of the norm polynomial of size $O_n(B)$,
while division by $p$ saves a factor of order $B/A$.

Before choosing this prime, consider $R<D^\epsilon$.
Substituting $A=B/R$
in \eqref{eq:AB} gives
\[
 B^{n-1}\ll_n R^mD^{-1/2},
\]
and consequently
\[
 H(f_1)\ll_n DB^n
 \ll_n D^{(n-2)/(2(n-1))}R^{nm/(n-1)}
 \ll_n D^e
\]
by \eqref{eq:epsilon}. This handles the second small-parameter
case.

We may therefore assume that $R\geq D^\epsilon$.
The scale and height are again admissible, since
\eqref{eq:AB} gives
\[
 R=\frac BA\ll_n\frac{D^{1/n}}X\ll_n D^{1/n},
 \qquad H(f_1)\ll_n DB^n\ll_n DX^{-n}\ll_n D.
\]
Lemma~\ref{lem:prime}, with $Y=R$, therefore supplies a
prime $R<p\leq2R$ avoiding the leading coefficient and
discriminant of $f_1$, together with a root modulo $p$.
For a root $0\leq r<p$, the polynomial
\[
 g(x)=\frac{f_1(px+r)}p
 =\frac Dp\N_{K/\Q}\bigl(p\alpha x+r\alpha-\beta\bigr)
\]
has integral coefficients and defines $K$ by
\eqref{eq:rescale}. Each conjugate of $p\alpha$ has
absolute value at most $pA$, and each conjugate of
$r\alpha-\beta$ has absolute value at most $pA+O_n(B)$.
The product bound therefore gives
\begin{equation}\label{eq:dual-saving}
 H(g)\ll_n\frac Dp(Ap+B)^n
 \ll_n DA B^{n-1}.
\end{equation}
Here $pA\asymp B$ and $1/p\asymp A/B$ explain the last
inequality. The calculation in \eqref{eq:dual-refined} now gives
$H(g)\ll_n D^aX^{-b}$.

Combining this with \eqref{eq:primal-qrh}, and applying
\eqref{eq:balance}, yields
\[
 \HH(K)\ll_{n,N}\min\{X^{n-1},D^aX^{-b}\}
 \leq D^{\rho(n,m)}.
\]
When $K$ is primitive, $m=1$, and the exponent is
$(n-1)/(2n)$.
\end{proof}

The proof works for any fixed $N\geq2$. The scales at which
we seek primes are fixed positive powers of $D$, and the
smaller scales have already been handled geometrically.
Thus $N$ affects the implied constant, but not the exponent.

\section{The optimal exponent for cubic fields}
\label{sec:cubic}

In this section, we give a direct proof, with an explicit constant, of the cubic lower-bound exponent established by Dubickas \cite[Theorem~1]{Dubickas}.

For an algebraic number $\alpha$ with primitive integral
minimal polynomial
\[
 f_\alpha(x)=T\prod_{j=1}^n(x-\alpha_j),\qquad T>0,
\]
its Mahler measure is
\[
 M(\alpha)=M(f_\alpha)
 =T\prod_{j=1}^n\max\{1,|\alpha_j|\}.
\]
For any polynomial $f$, define $M(f)$ by the same product
over its roots, including the absolute value of its leading
coefficient. Landau's inequality
\cite[Lemma~1.6.7]{BG} gives
\begin{equation}\label{eq:mahler-height}
 M(f)\leq\left(\sum_{i=0}^n|a_i|^2\right)^{1/2}
 \leq\sqrt{n+1}\,H(f).
\end{equation}
Jensen's formula and the inequality
between the geometric and arithmetic means give
\[
 M(f)^2
 =\exp\left(\frac1{2\pi}\int_0^{2\pi}
                 \log|f(e^{it})|^2\,dt\right)
 \leq\frac1{2\pi}\int_0^{2\pi}|f(e^{it})|^2\,dt
 =\sum_{i=0}^n|a_i|^2,
\]
where the last equality is Parseval's identity. The second
inequality in \eqref{eq:mahler-height} is immediate from
the definition of $H(f)$. Expanding the product over the
roots also gives the reverse comparison
\begin{equation}\label{eq:height-mahler}
 H(f)\leq2^nM(f);
\end{equation}
thus coefficient
height and Mahler measure differ by factors bounded in
terms of the degree alone.

We shall also use $M(\alpha^{-1})=M(\alpha)$ for
$\alpha\ne0$. To see this, the reciprocal polynomial
$x^nf_\alpha(1/x)$ is primitive and, up to sign, is the
minimal polynomial of $\alpha^{-1}$. Replacing each nonzero
root $\alpha_j$ by its reciprocal and the leading coefficient
by the constant coefficient preserves the defining product
for Mahler measure.

\begin{proposition}\label{prop:pure-cubic}
For every prime $p\equiv2\pmod9$, the field
$K_p=\Q(p^{1/3})$ satisfies
\[
 D_{K_p}=27p^2,
 \qquad
 \frac16D_{K_p}^{1/3}\leq\HH(K_p)\leq D_{K_p}^{1/3}.
\]
\end{proposition}

\begin{proof}
Put $\theta=p^{1/3}$. The polynomial $x^3-p$ is Eisenstein
at $p$, so $[K_p:\Q]=3$. The integral basis calculation
for pure cubic fields gives
\begin{equation}\label{eq:pure-basis}
 \mathcal O_{K_p}=\Z[\theta],\qquad D_{K_p}=27p^2;
\end{equation}
see \cite[Chapter~2, Exercise~41]{Marcus} or
\cite[Lemma~3]{Dubickas}. The condition $p\equiv2\pmod9$
ensures that $\Z[\theta]$ is the full ring of integers,
so its polynomial discriminant $-27p^2$ is also the field
discriminant.

Let $\alpha$ be any generator of $K_p$, and write
\[
 \alpha=b_0+b_1\theta+b_2\theta^2,\qquad b_i\in\Q.
\]
We may suppose that $b_2\ne0$, after replacing $\alpha$
by $\alpha^{-1}$ if necessary. Indeed, if $b_2=0$, then
$b_1\ne0$ because $\alpha$ is not rational, and
\[
 \alpha^{-1}
 =\frac{b_0^2-b_0b_1\theta+b_1^2\theta^2}
        {b_0^3+pb_1^3}
\]
has a nonzero coefficient of $\theta^2$. Its denominator
is $\N_{K_p/\Q}(\alpha)$, which is nonzero. Inversion
preserves the generated field and the Mahler measure, so a
lower bound after this replacement is also a lower bound
for the original generator.

Let $T>0$ be the leading coefficient of the primitive
integral minimal polynomial of the chosen $\alpha$. Write
this polynomial as $Tx^3+c_1x^2+c_2x+c_3$.
Multiplying its equation at $\alpha$ by $T^2$ shows that
$u=T\alpha$ satisfies
\[
 u^3+c_1u^2+Tc_2u+T^2c_3=0.
\]
Thus $T\alpha$ is an algebraic integer. By
\eqref{eq:pure-basis}, its coordinates $Tb_0,Tb_1,Tb_2$
in the basis $1,\theta,\theta^2$ are integers. In particular,
\[
 |Tb_2|\geq1.
\]
This step accounts for arbitrary denominators in a generator
and is the reason the lower bound applies to nonmonic
defining polynomials as well.

Let $\omega=e^{2\pi i/3}$, and write the conjugates as
\[
 \alpha_j=b_0+b_1\theta\omega^j
              +b_2\theta^2\omega^{2j},\qquad j=0,1,2.
\]
Using $1+\omega+\omega^2=0$, we isolate its quadratic
coordinate:
\[
 b_2\theta^2=\frac13\sum_{j=0}^2\omega^j\alpha_j.
\]
The triangle inequality therefore gives
$|b_2|p^{2/3}\leq\frac13\sum_j|\alpha_j|
\leq\max_j|\alpha_j|$. Hence
\[
 M(\alpha)
 =T\prod_{j=0}^2\max\{1,|\alpha_j|\}
 \geq T\max_j|\alpha_j|
 \geq T|b_2|p^{2/3}
 \geq p^{2/3}.
\]
Every integral defining cubic is an integer multiple of the
primitive minimal polynomial of its root, by Gauss's lemma.
Its height is at least the height of that primitive polynomial.
Equation~\eqref{eq:mahler-height}, with $n=3$, therefore gives
\[
 \HH(K_p)\geq\frac12p^{2/3}
 =\frac16D_{K_p}^{1/3}.
\]

For the upper bound, take $t=\lfloor p^{1/3}\rfloor$.
The polynomial
\[
 (x+t)^3-p=x^3+3tx^2+3t^2x+(t^3-p)
\]
defines $K_p$. Since
\[
 0\leq p-t^3
 =(p^{1/3}-t)(p^{2/3}+p^{1/3}t+t^2)
 \leq3p^{2/3},
\]
all its coefficients are at
most $3p^{2/3}=D_{K_p}^{1/3}$ in absolute value.
\end{proof}

Dirichlet's theorem \cite[Theorem~49]{Marcus} gives infinitely many primes
$p\equiv2\pmod9$. Their field discriminants are distinct,
so Proposition~\ref{prop:pure-cubic} supplies infinitely
many fields attaining the exponent $1/3$. This proves the
lower-bound assertion in Theorem~\ref{thm:cubic}.

\begin{proof}[Proof of Theorem~$\ref{thm:cubic}$]
Every cubic field is primitive, so Theorem~\ref{thm:qrh}
with $n=3$ gives $\HH(K)\ll_N D_K^{1/3}$ under the
stated zero-free hypothesis. Proposition~\ref{prop:pure-cubic}
gives the matching lower bound along an infinite family.
Since the discriminants in that family tend to infinity,
no smaller exponent can hold uniformly for all cubic fields.
\end{proof}

The same construction has a consequence in every degree
divisible by three. Passing to a fixed base field preserves the cubic
obstruction, provided we keep track of denominators over
that field. 

\begin{proposition}\label{prop:relative-cubic}
Fix a number field $F$ of degree $m\geq1$. For primes
$p\equiv2\pmod9$ unramified in $F$, put $K_p=F(p^{1/3})$.
Then
\[
 [K_p:\Q]=3m,\qquad \HH(K_p)\asymp_F D_{K_p}^{1/3}.
\]
In particular, for every positive integer $m$, there are
infinitely many fields of degree $n=3m$ with
$\HH(K)\asymp_n D_K^{1/3}$. If $m$ is odd, these fields
have no index-two subfield.
\end{proposition}

\begin{proof}
Put $\theta=p^{1/3}$ and $C_p=\Q(\theta)$. Since $p$ is
unramified in $F$, the polynomial $x^3-p$ is Eisenstein
at every prime of $F$ above $p$. Thus $[K_p:F]=3$, and
$F$ and $C_p$ are linearly disjoint.
Proposition~\ref{prop:pure-cubic} gives
$\mathcal O_{C_p}=\Z[\theta]$ and $D_{C_p}=27p^2$.
Set $R=\mathcal O_F[\theta]$ and $c=\gcd(D_F,27)$.
The integral-basis theorem for composita gives
\begin{equation}\label{eq:relative-basis}
 R\subseteq\mathcal O_{K_p}\subseteq c^{-1}R,
 \qquad |\Disc(R)|=D_F^3(27p^2)^m;
\end{equation}
see \cite[Chapter~2, Theorem~12, Corollary~1,
and Exercise~23(b)--(c)]{Marcus}. Here the common divisor of
$D_F$ and $D_{C_p}$ is $c$, since $p\nmid D_F$.
It follows that $[\mathcal O_{K_p}:R]\leq c^{3m}$, and
the discriminant-index formula yields
\begin{equation}\label{eq:relative-discriminant}
 D_{K_p}=\frac{D_F^3(27p^2)^m}{[\mathcal O_{K_p}:R]^2}
 \asymp_F p^{2m}.
\end{equation}
When $F$ is unramified at $3$, we have $c=1$, giving the
exact integral basis and discriminant used in the cubic case.

Let $\alpha$ be an arbitrary generator of $K_p$ over $\Q$,
and write
\[
 \alpha=b_0+b_1\theta+b_2\theta^2,\qquad b_i\in F.
\]
As in the cubic proof, either $\alpha$ or $\alpha^{-1}$
has $b_2\ne0$. Indeed, if $b_2=0$, then $b_1\ne0$
because $\alpha\notin F$, and the same reciprocal formula
applies with $b_0,b_1\in F$. We make this replacement if
necessary; it preserves $M(\alpha)$.

Let
\[
 g(x)=x^3+c_1x^2+c_2x+c_3\in F[x]
\]
be the monic minimal polynomial of $\alpha$ over $F$.
Define the integral ideal
\[
 \mathfrak a=\{u\in\mathcal O_F:
                  uc_i\in\mathcal O_F\text{ for }i=1,2,3\}.
\]
For each $u\in\mathfrak a$, the element $u\alpha$
satisfies
\[
 (u\alpha)^3+(uc_1)(u\alpha)^2
                  +(u^2c_2)(u\alpha)+u^3c_3=0.
\]
All coefficients are in $\mathcal O_F$, so
$\mathfrak a\alpha\subseteq\mathcal O_{K_p}$.
The inclusion in \eqref{eq:relative-basis} then gives
$c\mathfrak a\alpha\subseteq R$. Comparing the coefficients
of $1,\theta,\theta^2$, we obtain
$c\mathfrak a b_2\subseteq\mathcal O_F$. Taking the norm
of this nonzero integral ideal gives
\begin{equation}\label{eq:relative-denominator}
 \N_{F/\Q}(\mathfrak a)|\N_{F/\Q}(b_2)|\geq c^{-m}.
\end{equation}

The point is that $\N_{F/\Q}(\mathfrak a)$ is exactly
the leading coefficient $T>0$ of the primitive integral
minimal polynomial $f_\alpha$ over $\Q$. We verify this
standard consequence of Gauss's lemma for nonarchimedean
coefficient norms; see \cite[Lemma~1.6.3]{BG}.
Since $\alpha$ generates $K_p$ over $\Q$, we have
\[
 \frac{f_\alpha(x)}T
   =\prod_{\tau:F\hookrightarrow\overline\Q}g^\tau(x).
\]
For a prime ideal $\mathfrak p$ of $F$, put
\[
 a_{\mathfrak p}
 =\max\{0,-v_{\mathfrak p}(c_1),
              -v_{\mathfrak p}(c_2),-v_{\mathfrak p}(c_3)\},
\]
where $v_{\mathfrak p}$ takes integer values. Then
$\mathfrak a=\prod_{\mathfrak p}\mathfrak p^{a_{\mathfrak p}}$.
At a rational prime $q$, use the absolute value with
$|q|_q=q^{-1}$, and define the coefficient norm of a
polynomial to be the maximum of the $q$-adic absolute
values of its coefficients. An embedding
$F\hookrightarrow\overline{\Q}_q$ belonging to a prime
$\mathfrak p\mid q$ contributes coefficient norm
$q^{a_{\mathfrak p}/e(\mathfrak p/q)}$, and there are
$e(\mathfrak p/q)f(\mathfrak p/q)$ such embeddings.
Multiplicativity of the coefficient norm, together with
primitivity of $f_\alpha$, gives
\[
 v_q(T)=\sum_{\mathfrak p\mid q}
                   f(\mathfrak p/q)a_{\mathfrak p}.
\]
These are precisely the prime exponents in
$\N_{F/\Q}(\mathfrak a)$, proving the claimed equality.
Equation~\eqref{eq:relative-denominator} consequently reads
\begin{equation}\label{eq:relative-leading}
 T|\N_{F/\Q}(b_2)|\geq c^{-m}.
\end{equation}

Put $\omega=e^{2\pi i/3}$. For each embedding
$\tau:F\hookrightarrow\C$, the three
extensions to $K_p$ give the conjugates
\[
 \alpha_{\tau,j}
 =\tau(b_0)+\tau(b_1)\theta\omega^j
                +\tau(b_2)\theta^2\omega^{2j},
 \qquad j=0,1,2.
\]
The Fourier identity used above gives
\[
 \max_j|\alpha_{\tau,j}|
       \geq|\tau(b_2)|p^{2/3}.
\]
Multiplying over the $m$ embeddings of $F$ and using
\eqref{eq:relative-leading}, we obtain
\begin{align*}
 M(\alpha)
 &=T\prod_\tau\prod_{j=0}^2\max\{1,|\alpha_{\tau,j}|\}\\
 &\geq T|\N_{F/\Q}(b_2)|p^{2m/3}
 \geq c^{-m}p^{2m/3}.
\end{align*}
By \eqref{eq:mahler-height} and \eqref{eq:relative-discriminant},
this proves $\HH(K_p)\gg_F D_{K_p}^{1/3}$.

For the matching upper bound, fix an integral generator
$\eta$ of $F$ and set $t=\lfloor p^{1/3}\rfloor$.
For all sufficiently large $p$, the element
$\beta=\eta+\theta-t$ generates $K_p$ over $\Q$.
Its conjugates are
$\tau(\eta)+\theta\omega^j-t$. If two have the same
$j$, they are distinct because $\eta$ generates $F$.
If their $j$'s differ, the difference of their cubic terms
has size $\sqrt3\,p^{1/3}$, which eventually exceeds
every difference of the fixed conjugates of $\eta$.
Thus all $3m$ conjugates are distinct.

For each $\tau$, the conjugate with $j=0$ is bounded in
terms of $F$, since $0<\theta-t<1$, while the other two
are $O_F(p^{1/3})$. The element $\beta$ is integral, so
$M(\beta)\ll_F p^{2m/3}$. Equation~\eqref{eq:height-mahler}
gives $\HH(K_p)\ll_F D_{K_p}^{1/3}$.
The finitely many smaller permitted primes are absorbed
into the implied constant. There are infinitely many
permitted primes by Dirichlet's theorem. The prime $p$
ramifies in $K_p$, while all other ramified rational primes
divide $3D_F$; hence these fields are distinct. Choosing
$F$ once for each $m$ gives constants depending only on
$n=3m$. If $m$ is odd, the degree is odd, so there can be
no index-two subfield.
\end{proof}

\subsection*{Acknowledgments}

We are grateful to Akshay Venkatesh for many helpful conversations and, in particular, for explaining to us the results in~\cite{MV} and \cite{Venkatesh}.  The  authors were supported by the National Science Foundation (Grants DMS-1001828 and DMS-2303211). Much of this work was carried out while the authors were visiting the Lodha Mathematical Sciences Institute during December 2025; we thank LMSI for its hospitality.

\end{document}